\documentclass[12pt]{amsart}
\usepackage{amsfonts,amsmath,amsthm,amssymb}
\usepackage{latexsym}
\newtheorem{theorem}{Theorem}[section]
\newtheorem{corollary}[theorem]{Corollary}
\newtheorem{construction}[theorem]{Construction}
\newtheorem{lemma}[theorem]{Lemma}

\newtheorem{conjecture}[theorem]{Conjecture}
\newtheorem{proposition}[theorem]{Proposition}

\newtheorem{definition}{Definition}[section]

\theoremstyle{definition}

\theoremstyle{remark}
\newtheorem{rem}[theorem]{Remark}
\theoremstyle{remark}

\newcommand{\beql}[1]{\begin{equation}\label{#1}}
\newcommand{\eeq}{\end{equation}}
\newcommand \vs {{\bf s}}

\newcommand \vr {{\bf r}}
\newcommand \vq {{\bf q}}
\newcommand \vt {{\bf t}}
\newcommand \vu {{\bf u}}
\newcommand \vv {{\bf v}}
\newcommand \vw {{\bf w}}

\begin{document}

\title{Parametrizing complex Hadamard matrices}

\author{Ferenc Sz\"oll\H{o}si}

\address{Ferenc Sz\"oll\H{o}si: Budapest University of Technology and Economics (BUTE), H-1111, Egry J. u. 1, Budapest, Hungary}\email{szoferi@math.bme.hu}




\date{May, 2007.}

\begin{abstract}
The purpose of this paper is to introduce new parametric families
of complex Hadamard matrices in two different ways. First, we
prove that every real Hadamard matrix of order $N\geq 4$ admits an
affine orbit. This settles a recent open problem of Tadej and
\.Zyczkowski \cite{karol}, who asked whether a real Hadamard
matrix can be isolated among complex ones. In particular, we apply
our construction to the only (up to equivalence) real Hadamard
matrix of order $12$ and show that the arising affine family is
different from all previously known examples listed in
\cite{karol}. Second, we recall a well-known construction related
to real conference matrices, and show how to introduce an affine
parameter in the arising complex Hadamard matrices. This leads to
new parametric families of orders $10$ and $14$. An interesting
feature of both of our constructions is that the arising families
cannot be obtained via Di\c{t}\u{a}'s general method \cite{dita}.
Our results extend the recent catalogue of complex Hadamard
matrices \cite{karol}, and may lead to direct applications in
quantum-information theory.
\end{abstract}

 \maketitle



\section{Introduction}\label{sec:intro}
In the past few decades complex Hadamard matrices  were
extensively studied since it turned out that they are related to
many interesting combinatorial and important physical problems.
However, despite of many years of research only moderate results
are known, e.g. the problem of finding all complex Hadamard
matrices even of small orders is still open. The first significant
result is due to Haagerup \cite{haagerup}, who managed to classify
all complex Hadamard matrices up to order $5$ in $1997$. Only
partial results are known about matrices of order $6$. Besides
some affine families listed in \cite{karol}, all self-adjoint
(Hermitian) complex Hadamard matrices of order $6$ were classified by Beauchamp and
Nicoara \cite{nic}, and a symmetric non-affine family was found by
Matolcsi and Sz\"oll\H{o}si very recently \cite{sym}.

First, there was an interest in particular examples of
(permutation) inequivalent complex Hadamard matrices of low order.
However, due to a recent discovery of Di\c{t}\u{a} \cite{dita} the
situation has changed dramatically. His powerful method leads to
the construction of \emph{parametric} families of Hadamard
matrices in composite dimensions. This method was subsequently
rediscovered by Matolcsi, R\'effy and Sz\"oll\H{o}si \cite{MM} who
used a spectral set construction from \cite{KM}, and then used
another spectral set construction to obtain new families of
complex Hadamard matrices. An entirely different approach for
parametrization was described in the monument paper of Tadej and
\.Zyczkowski \cite{karol} who introduced the method of ``linear
variation of phases'', obtaining affine Hadamard families. They
successfully obtained all maximal affine Hadamard families
stemming from the Fourier matrices $F_N$ for $N\leq 16$. Thus, one
is interested in the inequivalent classes of {\emph{ parametric
families}} of Hadamard matrices nowadays.

The aim of this paper is to describe two \emph{general}
constructions which lead to new parametric families of complex
Hadamard matrices in certain dimensions; these matrices arise due
to a natural construction from real Hadamard and real conference
matrices. We prove that they are non-Di\c{t}\u{a}-type, which
subsequently leads to new results in the sense that they were not
included in the recent catalogue. The main point of this paper is
to show that these matrices always admit an affine orbit, thus we
can introduce new parametric families of complex Hadamard matrices
of order $10, 12$ and $14$. With the aid of our results we can
supplement the incomplete catalogue of complex Hadamard matrices
of small orders in \cite{karol}.

\section{Preliminaries}\label{prel}
First let us introduce some formal definitions  and recall
previous results from \cite{dita}, \cite{MM} and \cite{karol}.

\begin{definition}
An Hadamard matrix $H$ is a square complex matrix of order $N$
with $|H_{i,j}|=1$ for $i,j=1,2,\hdots, N$, satisfying $HH^\ast=NI$,
where $I$ is the identity matrix and $H^\ast$ denotes the
Hermitian transpose of $H$.
\end{definition}

\begin{definition}
A complex (real) Hadamard matrix $H$ of order $N$ is dephased
(normalized) if $H_{1,i}=H_{i,1}=1$ for every $i=1,2,\hdots,N$. In
a given dephased matrix $H$, the  lower right $(N-1)\times (N-1)$
submatrix is called the core of $H$.
\end{definition}

\begin{definition}\label{d}
Two Hadamard matrices, $H_1$ and $H_2$, are  equivalent if there
exist diagonal unitary matrices $D_1$ and $D_2$ and permutation
matrices $P_1$ and $P_2$ such that $H_1=D_1P_1H_2P_2D_2$.
\end{definition}

It is clear that every complex Hadamard matrix is equivalent to a
dephased one.

Next we recall Di\c{t}\u{a}'s general method of  constructing
complex Hadamard matrices (his subsequent results on families with
some free parameters follow easily from this formula as described
very well in his paper \cite{dita}).
\begin{construction}
Let $M$ be a complex Hadamard matrix of order $k$, and
$N_1,N_2,\hdots, N_{k}$ are  complex Hadamard matrices of order
$n$. Then \beql{ditaformula} K:=\left[
\begin{array}{cccc}
m_{11}N_{1}&\cdot&\cdot&m_{1k}N_k\\
\cdot&\cdot&\cdot&\cdot\\
\cdot&\cdot&\cdot&\cdot\\
m_{k1}N_1&\cdot&\cdot&m_{kk}N_k
\end{array}\right]
\eeq
is a complex Hadamard matrix of order $nk$.
\end{construction}

\begin{definition}
A complex Hadamard matrix $K$ is called Di\c{t}\u{a}-type if  it
is equivalent to a matrix arising from formula
\eqref{ditaformula}.
\end{definition}

\begin{definition}
A parametric family of complex Hadamard matrices is called affine if the phases of the entries are sums of a constant and a linear function of the parameters. A family is maximal affine, if it is not properly contained in any other affine family.
\end{definition}

\begin{rem}
When we say that $H$ admits an affine orbit, we mean that there exists an affine family stemming from a dephased form of $H$, consisting purely of dephased complex Hadamard matrices. Since the first row and column entries are fixed at some chosen values, the members of the family cannot be obtained by multiplication by unitary diagonal matrices.
\end{rem}

Several affine families are listed in \cite{karol}. For an example
of an affine family in this paper the reader might want to jump
ahead to formulas \eqref{h12}, \eqref{h7}, \eqref{r7}.

In general, deciding whether two Hadamard matrices are equivalent or not is a nontrivial task. However, recently Matolcsi et al. introduced a powerful method, which easily establishes if an Hadamard matrix is a \emph{Di\c{t}\u{a}-type} one. In fact, it turned out that it is worth investigating the corresponding log-Hadamard matrix. (A square matrix $L$ is log-Hadamard if the entrywise exponential matrix, $[e^{2\pi\textbf{i}L_{i,j}}]$, $L_{i,j} \in [0,1)$, is Hadamard). The following definition and Lemma \ref{RJ} summarize the corresponding results from \cite{MM}.

\begin{definition}\label{NNR}
Let $L$ be an $N\times N$ real matrix.  For an index set
$I=\{i_1,i_2,\dots, i_n\} \subset \{1,2,\dots, N\}$ two rows (or
columns) $\vs$ and $\vq$ are called {\it $I$-equivalent}, in
notation $\vs\sim_I\vq$, if the positive fractional part of the
entry-wise differences, $s_i-q_i \mod\ 1$, are the same for every
$i\in I$. Two rows (or columns) $\vs$ and $\vq$ are called
$(d)$-$n$-equivalent if there exist $n$-element disjoint sets of
indices $I_1, \dots , I_d$ such that $\vs\sim_{I_j}\vq$ for all
$j=1, \dots , d$.
\end{definition}

\begin{lemma}\label{RJ}
Permutation of rows and columns,  or adding a constant to a row or
a column does not change $(d)$-$n$-equivalence. \end{lemma}

By formula \eqref{ditaformula}, the structure of  an $N\times N$ Di\c{t}\u{a}-type matrix $L$ (where $N=nk$) implies for the corresponding log-Hadamard matrix $\log L$ that there exists a partition of indices into $n$-element sets $I_1, \dots , I_k$ and $k$-tuples of rows $R_j=\{\vr^j_1, \dots, \vr^j_k\}$ ($j=1, \dots, n$) such that any two rows in a fixed $k$-tuple are equivalent with respect to any of the $I_m$'s. Naturally, the same holds for the transpose of a Di\c{t}\u{a}-type matrix, with the role of rows and columns interchanged.

The following observation is a trivial consequence  of their result:

\begin{lemma}\label{O}
Let $H$ be a dephased complex Hadamard matrix of order $N$, and
suppose that $H_{i,j}\neq 1$ for every $1<i,j\leq N$, i.e.
there is no $1$ in the core of $H$. Then $H$ is not of
Di\c{t}\u{a}-type.
\end{lemma}

\begin{proof}
We argue by contradiction. Assume that $H$ is Di\c{t}\u{a}-type.
Using the notations of the previous paragraph we can arrange
(after relabelling the index sets if necessary) that
$\{1\}\subseteq I_1$ and (after permuting the columns of $H$ if
necessary) that $\{1,2\}\subseteq I_1$. There must be a row $\vr$
of $\log H$ which is $I_1$-equivalent to the first row. However,
as all entries in the first row and first column are $0$'s, this
would imply that $\vr$ contains a 0 in its second coordinate, a
contradiction.
\end{proof}

\section{Constructing complex Hadamard matrices from real ones}
In this section we investigate the structure of real Hadamard
matrices.  First we prove that they cannot be obtained using
Di\c{t}\u{a}'s method in certain dimensions. Next we introduce a
somewhat natural construction for obtaining new, parametrized
complex Hadamard matrices from real ones. In fact, it was asked in
\cite{karol} whether all real Hadamard matrices of order $N\geq4$
can be parametrized and, by Theorem \ref{b}, we answer this
question in the positive. Before doing so we first recall a
folklore

\begin{lemma}\label{L}
Let $p\geq 3$ be an arbitrary odd number.  Suppose that the first
four rows of a real $\{-1,1\}$ matrix of order $4p$ have the
following form (note that every real Hadamard matrix is easily
seen to be equivalent to one having exactly the same first {\it
three} rows as the matrix below): \beql{can}
\begin{array}{rrrr}
(\vs) \\
(\vt) \\
(\vu) \\
(\vv) \\
\end{array}
\left[\begin{array}{cccc}
1^p & 1^p & 1^p & 1^p\\
1^p & 1^p & (-1)^p & (-1)^p\\
1^p & (-1)^p & 1^p & (-1)^p\\
1^p & (-1)^p & (-1)^p & 1^p\\
\end{array}\right],
\eeq where $1^p$ means $p$ one's in a row.  Then this matrix
cannot be extended with a further $\{1,-1\}$ row being orthogonal to all previous ones.
\end{lemma}
\begin{proof}
Suppose, to the contrary, that \eqref{can} can be extended by a
further row $\vw$. Let us denote by $a,b,c$ and $d$ the number of
1's in $\vw$ in the first-, second-, third- and fourth quarter,
i.e.
$\vw=\left(1^a,(-1)^{p-a},1^b,(-1)^{p-b},1^c,(-1)^{p-c},1^d,(-1)^{p-d}\right),
0\leq a,b,c,d\leq p$. Since $\vw$ is orthogonal to all of the rows
$\vs, \vt, \vu$ and $\vv$, we get the following four equations by
straightforward computation

\beql{12}
a-(p-a)+b-(p-b)+c-(p-c)+d-(p-d)=0
\eeq

\beql{13}
a-(p-a)+b-(p-b)-c+(p-c)-d+(p-d)=0
\eeq

\beql{14}
a-(p-a)-b+(p-b)+c-(p-c)-d+(p-d)=0
\eeq

\beql{15}
a-(p-a)-b+(p-b)-c+(p-c)+d-(p-d)=0
\eeq

By simple algebra  one can check that the solution to equations
\eqref{12}-\eqref{15} is $a=b=c=d=\frac{p}{2}$ and, since $p$ is
odd by assumption, this is a contradiction.
\end{proof}

Now we are ready to state our first

\begin{theorem}\label{nRJ}
Let $p$ be an odd prime and suppose that $H_{4p}$ is a \emph{real} Hadamard matrix of order $4p$. Then $H_{4p}$ is not of Di\c{t}\u{a}-type.
\end{theorem}

\begin{proof}
We will use the notation of the paragraph following Lemma \ref{RJ}
with the exception that instead of taking $\log H$ we apply the
notion of $I$-equivalence to the rows of $H$ itself in a natural
way.

Assume, to the  contrary, that $H_{4p}$ is of Di\c{t}\u{a}-type.
In this case the only possible values for $n$ are $2,4,p$ and $2p$
(with $k$ being $2p,p,4$ and $2$ respectively). Suppose that
$H_{4p}$ is dephased, and let us again denote the rows of \eqref{can} by $\vs,\vt,\vu$ and $\vv$ respectively. There are four cases to consider according to the choices of $n$ and $k$:

\textbf{CASE 1}  Assume $n=2p, k=2$. In this case there should be
a partition of indices to $2p$-element sets $I_1, I_2$ such that
in $H_{4p}$ $2p$ pairs of rows are equivalent with respect to
$I_1$ and $I_2$.  After permutation of rows and columns it is
trivial to achieve that the first three rows of $H_{4p}$ are $\vs,
\vt$ and $\vu$, respectively, and $\vs$ and $\vt$ form a pair.
(First we permute the rows so that the companion of row $1$
becomes the second row and then we permute the columns so that the
position of $1$'s and $-1$'s is exactly as in \eqref{can}.) Then
$I_1=\{1,2,\hdots,2p\}$ and $I_2=\{2p+1,2p+2,\hdots, 4p\}$. Now
consider $\vu$. If it formed a pair, then its companion's first
$2p$ entries would have to be exactly the same as those in $\vu$.
However, by orthogonality, the last $2p$ entries in $\vu$ and its
companion must be opposite. Thus the companion of $\vu$ must be
exactly $\vv$, which is a contradiction since there is no such a
row in $H_{4p}$ due to Lemma \ref{L} (by our assumptions, of
course, $H$ has at least $12$ rows).

\textbf{CASE 2}  Now assume $n=p,k=4$. In this case the partitions
of indices are $p$-element sets $I_1,I_2,I_3$ and $I_4$, such that
in $H_{4p}$ there exists $p$ $4$-tuples of rows, such that any two
rows in a fixed $4$-tuple are equivalent with respect to them. We
can suppose that $I_1=\{1,2,\hdots, p\},
I_2=\{p+1,p+2,\hdots,2p\},I_3=\{2p+1,2p+2,\hdots,3p\}$ and
$I_4=\{3p+1,3p+2,\hdots, 4p\}$. Now observe, since $\vs$ contains
only $1$'s, any row equivalent to it with respect to $I_1,I_2,I_3$
and $I_4$ must be one of $\vt,\vu$ or $\vv$. However, we need
three rows being equivalent to $\vs$, thus we need all four rows of
\eqref{can}, which is a contradiction again.

\textbf{CASE 3}  Now assume $n=4, k=p$. In this case the
partitions of indices are $4$-element sets $I_1,I_2,\hdots ,I_p$
such that in $H_{4p}$ there exist $4$ disjoint $p$-tuples of rows
such that any two rows in a fixed $p$-tuple are equivalent with
respect to them. Again, we would like to find a companion to
$\vs$. Observe that since every row (different from $\vs$)
contains $2p$ $1$'s and $2p$ $(-1)$'s it is impossible to split
their entries into odd ($p$) number of disjoint sets containing
exactly the same values. Hence we cannot choose a companion to
$\vs$, equivalent to it with respect to the index sets.

\textbf{CASE 4}  Finally assume that $n=2, k=2p$. Again (by
permuting the columns of $H_{4p}$ if necessary), we can suppose
that $I_1=\{1,2\},I_2=\{3,4\}$. Since $H_{4p}$ is a real Hadamard
matrix, we can suppose that (after permuting some rows if
necessary) its first three columns are exactly the same as the
transpose of the first three row of \eqref{can}. Now observe that
the first $2p$ and the second $2p$ rows have to belong to a common
tuple. To preserve equivalence with respect to $I_2$, one can see
that the fourth column of $H_{4p}$ has to be exactly the transpose
of the fourth row of the matrix in \eqref{can}. And this is a
contradiction again.
\end{proof}

\begin{corollary}\label{imp}
$H_{12}$ is not of Di\c{t}\u{a}-type.
\end{corollary}

Now we turn to the parametrization of real Hadamard matrices. It
is well known that $H_4$ admits a $1$-parameter orbit. In
\cite{karol} a $5$-parameter, while in \cite{MM} a $4$-parameter
maximal affine orbit was constructed for $H_{8}$ (these orbits are
essentially different, but they intersect each other at $H_8$). In
general it is not clear how to introduce affine parameters to an
\emph{arbitrary} complex Hadamard matrix. The authors of
\cite{karol} admit that the ``linear variation of phases'' method
becomes a serious combinatorial problem already for $N = 12$, so
it cannot effectively be used for higher order matrices. Now we
introduce a general method for parametrization which always works
for real matrices and, in some cases, for complex matrices too.
The main observation is contained in the following

\begin{lemma}\label{TT}
Let $H$ be an arbitrary dephased complex Hadamard matrix of order
$N\geq 4$.  Suppose that $H$ has a pair of columns, say $\vu$ and
$\vv$, with the following property: $u_i=v_i$ or $u_i+v_i=0$ holds
for every $i=1,2,\hdots, N$. Then $H$ admits an affine orbit.
\end{lemma}

\begin{proof}
Consider $H$ satisfying the conditions of  Lemma \ref{TT}, and
take every pair of coordinates $(u_i,v_i)$ for which $u_i+v_i=0$
holds. Multiply these elements by $e^{\textbf{i}t}$, i.e. modify
$(u_i,v_i)$ to $(u_ie^{\textbf{i}t},v_ie^{\textbf{i}t})$. Now we
proceed to show that the arising parametric matrix $H^{(1)}(t)$ is
Hadamard. To do this let us consider a pair of rows in
$H^{(1)}(t)$. It is easy to see that after taking the inner
product of these rows, the parameter (if it existed in at least
one of them) vanishes, therefore $H^{(1)}(t)$ is Hadamard
independently of the exact value of $t$. Finally, if $H^{(1)}(t)$
is not dephased (i.e. we have chosen the first column of $H$ to be
either $\vu$ or $\vv$), one should multiply some rows by
$e^{-\textbf{i}t}$ to get a dephased matrix, and it is clear that
$t$ will not vanish whenever $N\geq 4$.
\end{proof}

With the aid of Lemma \ref{TT} we can prove  the main theorem of
this section. We prove that there is no isolated matrix among real
Hadamard matrices except for orders $1$ and $2$ (the cases $N=4$
and $N=8$ were mentioned in the paragraph preceding Lemma
\ref{TT}) .

\begin{theorem}\label{b}
Let $H$ be a real Hadamard matrix of order  $N\geq 12$. Then $H$
admits an $\left(\frac{N}{2}+1\right)$-parameter affine orbit.
\end{theorem}

\begin{proof}
Let $N\geq 12$, and let us take an arbitrary  dephased real
Hadamard matrix of order $N$, say $H$. It is clear that when
considering any two columns of $H$, there will be exactly
$\frac{N}{2}$ rows, where the entries of these columns differ, and
another $\frac{N}{2}$ rows, where the entries of these columns are
the same, so the conditions of Lemma \ref{TT} hold. Now we apply
the construction described in the proof of Lemma \ref{TT}.

Clearly, we can further assume, that $H$  has the following
``canonical'' form: $H_{2,1}=H_{2,2}=\hdots=H_{2,N/2}=1$ so
$H_{2,N/2+1}=H_{2,N/2+2}=\hdots=H_{2,N}=-1$ and
$H_{3,3}=H_{3,4}=1$ and $H_{3,2}=H_{3,N-1}=H_{3,N}=-1$. Consider
the following set containing pairs of indices: $T=\{(2i-1,2i) :
i=1,2,\hdots,\frac{N}{2}\}$. Every element of $T$ represents a
pair of columns in $H$. Now the construction is the following: for
every $i=1,2,\hdots, \frac{N}{2}$ take the respective element of
$T$, and consider the \emph{rows} of the corresponding pair of
columns. If the entries in a row are different then multiply them
by $e^{\textbf{i}x_i}$ (again: there are exactly $\frac{N}{2}$
such rows). This yields an $\frac{N}{2}$-parameter family,
stemming from $H$. However, it is not dephased, so one has to
multiply some rows by $e^{-\textbf{i}x_1}$ to get a dephased
Hadamard matrix. Since $H_{3,3}=H_{3,4}$ we can see that these
entries, after parametrization and dephasing the matrix, depend
only on $x_1$, so $x_1, x_2,\hdots,x_{N/2}$ are independent
parameters in the dephased matrix. For convenience, we can
substitute $x_1$ by $-x_1$. Now taking a look at the first two
\emph{rows} of $H$ (which are still independent, after
parametrization, of any of the $x_i$'s) one can multiply the last
(differing) $\frac{N}{2}$ entries of these by
$e^{-\textbf{i}x_{N/2+1}}$, the arising matrix thus being still
Hadamard. Again, it is not dephased, but observe that after
dephasing the matrix (multiplying the last $\frac{N}{2}$ columns
by $e^{\textbf{i}x_{N/2+1}}$), since $H_{3,N-1}=H_{3,N}$ these
entries after parametrization depend only on $x_1$ and on
$x_{N/2+1}$. Note that this last operation left unchanged both the
parametrized $H_{3,3}$ and $H_{3,4}$ which still depend only on
$x_1$. This completes the proof.
\end{proof}
\begin{rem}
The same construction also works when  we replace ``rows'' by ``columns'' and vice versa.
\end{rem}

\begin{rem}
It is easy to see (by taking the inner  product of $\vu$ and
$\vv$) that Lemma \ref{TT} can only be applied in even orders.
However, the conditions of this lemma hold for many non-real
Hadamard matrices, too. For example, the Fourier matrix $F_N$ in
even orders \emph{has} two columns in which the entries are either the same or of opposite sign. Other examples are the matrices $S_8, S_{12}$ and $S_{16}$ in \cite{MM} which also satisfy the conditions of Lemma \ref{TT}. Thus, this lemma can be used for parametrizing a wide class of complex Hadamard matrices.
\end{rem}

Now we give an example. The following  matrix is the only real
Hadamard matrix of order $12$ (up to equivalence). We note that it
can be constructed from a skew-symmetric conference matrix (see
section \ref{conf5}).

\beql{h12}
H_{12}=\left[\begin{array}{rrrrrrrrrrrr}
 1 & 1 & 1 & 1 & 1 & 1 & 1 & 1 & 1 & 1 & 1 & 1 \\
 1 & 1 & 1 & 1 & 1 & 1 & -1 & -1 & -1 & -1 & -1 & -1 \\
 1 & -1 & 1 & 1 & -1 & -1 & 1 & 1 & -1 & 1 & -1 & -1 \\
 1 & -1 & 1 & -1 & 1 & -1 & -1 & 1 & 1 & -1 & 1 & -1 \\
 1 & -1 & -1 & 1 & -1 & 1 & -1 & 1 & 1 & -1 & -1 & 1 \\
 1 & 1 & 1 & -1 & -1 & -1 & 1 & -1 & 1 & -1 & -1 & 1 \\
 1 & 1 & -1 & -1 & 1 & -1 & -1 & 1 & -1 & 1 & -1 & 1 \\
 1 & 1 & -1 & 1 & -1 & -1 & -1 & -1 & 1 & 1 & 1 & -1 \\
 1 & -1 & 1 & -1 & -1 & 1 & -1 & -1 & -1 & 1 & 1 & 1 \\
 1 & 1 & -1 & -1 & -1 & 1 & 1 & 1 & -1 & -1 & 1 & -1 \\
 1 & -1 & -1 & -1 & 1 & 1 & 1 & -1 & 1 & 1 & -1 & -1 \\
 1 & -1 & -1 & 1 & 1 & -1 & 1 & -1 & -1 & -1 & 1 & 1\\
\end{array}\right]
\eeq

By Theorem \ref{b} we can easily  construct a $7$-parameter family
stemming from $H_{12}$. The notations here are exactly the same as
in \cite{karol} and \cite{MM}. We denote by $\circ$ the Hadamard
product of two matrices (i.e. $[H_1\circ
H_2]_{i,j}=[H_1]_{i,j}\cdot[H_2]_{i,j}$), while the symbol
$\textbf{EXP}$ stands for the entrywise exponential operation (i.e.
$[\textbf{EXP}H]_{i,j}=exp(H_{i,j})$).

\beql{h7}
H^{(7)}_{12}(a,b,c,d,e,f,g)=H_{12}\circ\textbf{EXP}\left(\textbf{i}\cdot
R_{H_{12}^{(7)}}(a,b,c,d,e,f,g)\right) \eeq

where

\beql{r7} R_{H^{(7)}_{12}}(a,b,c,d,e,f,g)= \eeq

\scriptsize
\[\left[\begin{array}{cccccccccccc}
\bullet & \bullet & \bullet & \bullet & \bullet & \bullet & \bullet & \bullet & \bullet & \bullet & \bullet & \bullet\\
 \bullet & \bullet & \bullet & \bullet & \bullet & \bullet & \bullet & \bullet & \bullet & \bullet & \bullet & \bullet\\
 \bullet & \bullet & a & a & a & a & a+ g & a+ g & a+ e+ g & a+ e+ g & a+ g & a+ g \\
 \bullet & \bullet & a+ b & a+ b & a+ c & a+ c & a+ d+ g & a+ d+ g & a+ e+ g & a+ e+ g & a+ f+ g & a+ f+ g \\
 \bullet & \bullet & a+ b & a+ b & a+ c & a+ c & a+ d+ g & a+ d+ g & a+ e+ g & a+ e+ g & a+ f+ g & a+ f+ g \\
 \bullet & \bullet &  b &  b & \bullet & \bullet &  d+ g &  d+ g &  e+ g &  e+ g &  f+ g &  f+ g \\
 \bullet & \bullet & \bullet & \bullet &  c &  c &  d+ g &  d+ g &  e+ g &  e+ g &  f+ g &  f+ g \\
 \bullet & \bullet &  b &  b & \bullet & \bullet &  g &  g &  g &  g &  f+ g &  f+ g \\
 \bullet & \bullet & a+ b & a+ b & a+ c & a+ c & a+ g & a+ g & a+ e+ g & a+ e+ g & a+ g & a+ g \\
 \bullet & \bullet & \bullet & \bullet &  c &  c &  g &  g &  g &  g &  f+ g &  f+ g \\
 \bullet & \bullet & a & a & a & a & a+ d+ g & a+ d+ g & a+ g & a+ g & a+ g & a+ g \\
 \bullet & \bullet & a+ b & a+ b & a+ c & a+ c & a+ d+ g & a+ d+ g & a+ g & a+ g & a+ g & a+ g
\end{array}\right]\]

\normalsize According to Corollary \ref{imp},  $H_{12}$ is not of
Di\c{t}\u{a}-type, so it admits only non-Di\c{t}\u{a}-type
matrices in a small neighbourhood of it, since the the set of
Di\c{t}\u{a}-type matrices is closed as shown in the following

\begin{proposition}\label{NON}
The set of all $N\times N$ Di\c{t}\u{a}-type matrices is closed in
the space of all $N\times N$ matrices.
\end{proposition}

\begin{proof}
Let $T_l\to T$ be a convergent sequence of Di\c{t}\u{a}-type
matrices. We need to show that $T$ is also Di\c{t}\u{a}-type.

By definition there exist permutation matrices $P^{(l)}_1,
P^{(l)}_2$ and diagonal unitary matrices $D^{(l)}_1, D^{(l)}_2$
such that $P^{(l)}_1D^{(l)}_1T_lD^{(l)}_2P^{(l)}_2=K_l$, where
$K_l$ arises in formula \eqref{ditaformula}. Each $K_l$ can be
characterized by the values of $k$, $m_{11}, \dots, m_{kk}$, and
the matrices $N_1, \dots, N_k$ in \eqref{ditaformula} (each
depending on $l$, of course, which we left out to simplify the
notation). Since the number of possible permutation matrices and
the number of possible choices for $k$ is finite, and all other
parameters such as $D^{(l)}_1,$ $D^{(l)}_2,$ $m_{ij},$ $N_i$ take
values in compact spaces, there exists a subsequence $l_h$ along
which the permutation matrices and the value of $k$ are constant
and all other parameters converge, i.e. $D^{(l_h)}_1\to D_1$,
$D^{(l_h)}_2\to D_2$, $m^{l_h}_{ij}\to r_{ij}$, $N^{l_h}_i\to
Q_i$. By taking the limit it is clear that $T$ is equivalent to
the Di\c{t}\u{a}-type matrix $K$ characterized by the values $k$,
$r_{11}, \dots, r_{kk}$, and the matrices $Q_1, \dots, Q_k$ in
\eqref{ditaformula}.
\end{proof}

As a consequence, we have

\begin{corollary}
The family $H^{(7)}_{12}(a,b,c,d,e,f,g)$ contains only
non-Di\c{t}\u{a}-type  matrices in a small neighbourhood around $H_{12}$.
\end{corollary}

Now we show that $H_{12}$ is inequivalent  to any of the order
$12$ matrices appearing in \cite{karol} and \cite{MM}. First we
recall a result from Haagerup, who introduced the following set
$\Lambda_H=\{h_{ij}h_{kl}\overline{h}_{kj}\overline{h}_{il} :
(i,j,k,l) \in \{1, \hdots, N\}^{\times 4}\}$ for $H$ of order $N$.
In \cite{haagerup} he claims that this set is invariant under the
equivalence preserving operations, see Definition \ref{d}.

\begin{lemma}
Two complex Hadamard matrices, say $H_1$ and $H_2$, are
inequivalent, if they have different $\Lambda_H$-sets.
\end{lemma}

Now we are ready  to prove\footnote{The author is grateful to M.
Matolcsi who suggested the proof of Lemma \ref{g}.} the following

\begin{lemma}\label{g}
$H_{12}$ is inequivalent  to any of the $12\times 12$ matrices
listed in \cite{karol} and \cite{MM}.
\end{lemma}

\begin{proof}
The proof relies  on the Haagerup condition. First observe that
$\Lambda_{H_{12}}=\{1,-1\}$. Now consider the seven families of
order $12$ in \cite{karol} stemming from $F_{12}$, and notice that
$e^{2\pi \textbf{i}/8}\in \Lambda_{F_{12}}$ for any matrix of any
of these families stemming from $F_{12}$, independently of the
values of the parameters. Secondly, observe that
$e^{2\pi\textbf{i}/3}\in \Lambda_{S_{12}}$ for any matrix stemming
from $S_{12}$ in \cite{MM}, again independently of the actual
values of the parameters. These observations can be easily
verified by taking $h_{11}=\overline{h}_{1j}=\overline{h}_{i1}=1$,
and taking an appropriate element $h_{ij}$ for every matrix in the
families stemming from $F_{12}$ and from $S_{12}$. Since $H_{12}$
and matrices from these families possess different $\Lambda$-sets,
they cannot be equivalent. There are several other families of
order $12$ listed in \cite{karol}, however those families were
obtained by Di\c{t}\u{a}'s construction (and thus consist purely
of Di\c{t}\u{a}-type matrices), therefore they cannot contain a
non-Di\c{t}\u{a}-type matrix such as $H_{12}$. This completes the
proof.
\end{proof}

\begin{proposition}
The family $H^{(7)}_{12}(a,b,c,d,e,f,g)$ is  locally inequivalent
to the families presented in \cite{karol} and \cite{MM}.$\hfill
\square$
\end{proposition}

\begin{proof}
This clearly  follows from the fact, that the invariant set
$\Lambda$ changes continuously. If we change some entries in
$H_{12}$ from $\pm 1$ to $e^{\textbf{i}t}$ with
$0<\left|t\right|<\varepsilon$ or
$0<\left|t-\pi\right|<\varepsilon$ (for $\varepsilon$ being small)
then neither $e^{2\pi \textbf{i}/8}$ nor $e^{2\pi\textbf{i}/3}$
will arise in the $\Lambda$-set of the modified matrix. Finally,
by Proposition \ref{NON}, it is clear that we can choose
$\varepsilon$ small enough to obtain non-Di\c{t}\u{a}-type
matrices only.
\end{proof}

Finally, we consider dimension $16$.  The situation here is more
complicated since there are $5$ inequivalent real Hadamard
matrices of that order. Therefore, with the aid of our
construction (described in the proof of Theorem \ref{b}) we can
obtain $5$, locally inequivalent, parametrized families of complex
Hadamard matrices. The fact that parametric families stemming from
inequivalent Hadamard matrices are locally inequivalent can be
proved by the same argument as in Proposition \ref{NON}.

It is known that the orbit of the Fourier matrix $F_{16}$ passes
through one of the $5$ inequivalent real Hadamard matrices, namely the
matrix $F_2\otimes F_2\otimes F_2\otimes F_2$. Unfortunately we do
not know how the other $4$ real Hadamard matrices are related to
$F_{16}$ or to the recently constructed ``spectral set'' matrix
$S_{16}$ in \cite{MM}. However, as we mentioned before, $H_8$ can
be parametrized in at least two essentially different ways, and
that is exactly why we conjecture that the parametrized complex
Hadamard matrices constructed by Theorem \ref{b} are, at least
locally, new.

\section{Constructing complex Hadamard matrices from conference matrices}\label{conf5}
The aim of this section is to describe another general method for
constructing  parametrized complex Hadamard matrices. First we
recall a well-known and widely studied class of matrices:

\begin{definition}
A conference matrix of order $N$ is a square $N\times N$  matrix
$C$, satisfying $CC^T=C^TC=(N-1)I$, $C_{ii}=0, i=1,2,\hdots, N$ and
$C_{ij}\in\{-1,1\}$ for $i\neq j$.
\end{definition}

It is easy to see that for a given conference matrix $C$ either
multiplying any row or column by $-1$, or permuting the rows and
columns of $C$ with the same permutation matrix $P$ (i.e.
considering $PCP^T$ instead of $C$) we get a conference matrix
again. Conference matrices related in these two ways are called
equivalent. It is a well-known fact that real conference matrices
lead to an obvious construction of Hadamard matrices. Whenever $C$
is a real \emph{symmetric} conference matrix, then
`$H=I+\textbf{i}C$' is a complex Hadamard matrix. (For
skew-symmetric conference matrices the formula `$H=I-C$' is used).
In the rest of this paper we will refer to the `$H=I+\textbf{i}C$'
construction as the \emph{conference matrix construction}. It is
clear that equivalent conference matrices give rise to equivalent
Hadamard matrices. For a survey on conference matrices see e.g.
\cite{ray} or \cite{Cr}. There are infinitely many orders for
which a symmetric conference matrix exists, however it is still an
open problem to give a full characterization of them; it is well
known that the order of a conference matrix must be even, moreover
the order of a \emph{symmetric} conference matrix must be $N=4k+2$
for some nonnegative integer $k$. However this condition is not
sufficient due to a negative result proved by Raghavarao in
\cite{R}. In particular, if $N$ is the order of a symmetric
conference matrix, then $N-1$ must be the sum of two squares.
For a more or less up-to-date list of the orders of the known
conference matrices see the last sections of \cite{JS}.

Next we prove a general method for introducing an affine
parameter to every complex Hadamard matrix arising from the
conference matrix construction. We denote this class of complex
Hadamard matrices by $D$, as $D_6$ in \cite{karol} is exactly a
matrix arising from a symmetric conference matrix of order $6$.
The following statements are analogous to Theorem \ref{nRJ} and
Theorem \ref{b}.

\begin{theorem}\label{CO2}
Complex Hadamard matrices arising from the conference matrix
construction are not of Di\c{t}\u{a}-type.
\end{theorem}

\begin{proof}
After dephasing $H=I+\textbf{i}C$, the core of the resulting matrix
will contain $-1$'s in the main diagonal and $\pm \textbf{i}$'s
otherwise, therefore the statement follows from Lemma \ref{O}.
\end{proof}

\begin{theorem}\label{CO}
Every complex Hadamard matrix $D_N$ arising from the conference
matrix construction admits an affine orbit, i.e. there exists an
affine family of complex Hadamard matrices of at least one
parameter which contains $D_N$.
\end{theorem}

\begin{proof}
The proof is completely elementary, but requires many  cases to
consider. Let $D_N$ be any matrix arising from the conference
matrix construction, of order $N$. Further, we can arrange that it
be both symmetric and dephased (of course, after parametrization,
$D_N$ can be transformed back to the original form
$I+\textbf{i}C$, and this transformation clearly does not affect
the presence of parameters). In \cite{dita} and \cite{karol}
$D_6^{(1)}(t)$ appeared, as a parametric family of order $6$, so
we restrict our attention to the next order $N=4k+2$, and we
suppose that $N\geq 10$. We show that one parameter can be
introduced independently of what a conference matrix $C$ was used
to construct $D_N$. Indeed, consider its second (\vu) and third
(\vv) rows. Because $D_N$ is Hadamard, there are exactly
$\frac{N-2}{2}$ places where the entries of $\vu$ and $\vv$ differ
only by a sign. Multiply these entries by $e^{\textbf{i}t}$. Now
consider the second and the third \emph{column} of $D_N$, and
multiply those entries by $e^{-\textbf{i}t}$ which differ by a
sign row-wise. We prove that the obtained $1$-parameter matrix
$D_N^{(1)}(t)$  will still be Hadamard. We show that the modified
rows of $D_N$ are orthogonal to each and every other row of
$D_N^{(1)}(t)$ independently of $t$. There are many trivial cases,
but there are two which require some extra considerations:

\textbf{CASE 1:} We proceed to show that both $\vu$ and $\vv$ are
orthogonal  to any unchanged row. After permuting the rows and the
columns of $D_N^{(1)}(t)$, we can suppose that it has the
following (symmetric) form as beneath; it is also clear, that (by
taking the Hermitian transpose of $D_N^{(1)}(t)$ if it is
necessary and, again, permuting) the imaginary elements in the
upper left $3\times 3$ submatrix are $\textbf{i}$'s. Now consider
any unchanged row, other than the first row of $D_N^{(1)}(t)$; its
first three elements could be either $(1,\textbf{i},\textbf{i})$
or $(1,-\textbf{i},-\textbf{i})$ respectively. We consider the
first case, the other could be treated exactly in the same way.
Below in the figure one can see a sketch of $D_N^{(1)}(t)$.

\scriptsize

\[
\begin{array}{r}
(\vu)\\
(\vv)\\
\\
\\
\\
\\
\\
\\
\\
\\
\\
\\
\\
\\
\\
\\
\\
\end{array}
\left[\begin{array}{r|rr|r|rrr|rrr|rrr|rrr}
1 & 1 & 1 & 1 & 1 & \hdots & 1 & 1 & \hdots & 1 & 1 & \hdots & 1 & 1 & \hdots & 1\\
\hline
1 & -1 & \textbf{i} & \textbf{i} & \textbf{i} & \hdots & \textbf{i} & \textbf{i}e^{\textbf{i}t} & \hdots & \textbf{i}e^{\textbf{i}t} & -\textbf{i}e^{\textbf{i}t} & \hdots & -\textbf{i}e^{\textbf{i}t} & -\textbf{i} & \hdots & -\textbf{i}\\
1 & \textbf{i} & -1 & \textbf{i} & \textbf{i} & \hdots & \textbf{i} & -\textbf{i}e^{\textbf{i}t} & \hdots & -\textbf{i}e^{\textbf{i}t} & \textbf{i}e^{\textbf{i}t} & \hdots & \textbf{i}e^{\textbf{i}t} & -\textbf{i} & \hdots & -\textbf{i} \\
\hline
1 & \textbf{i} & \textbf{i} & -1 & $a$ & & $b$ & $c$ & & $d$ & $e$ & & $f$ & $g$ & & $h$\\
\hline
1 & \textbf{i} & \textbf{i} & & -1 & & & & & & & & & & &\\
\vdots & \vdots & \vdots & & & \ddots& & & & & & & & & & \\
1 & \textbf{i} & \textbf{i} & & & & -1 & & & & & & & & & \\
\hline
1 & \textbf{i}e^{-\textbf{i}t} & -\textbf{i}e^{-\textbf{i}t} & & & & & -1 & & & & & & & & \\
\vdots & \vdots & \vdots & & & & & & \ddots & & & & & & & \\
1 & \textbf{i}e^{-\textbf{i}t} & -\textbf{i}e^{-\textbf{i}t} & & & & & & & -1 & & & & & & \\
\hline
1 & -\textbf{i}e^{-\textbf{i}t} & \textbf{i}e^{-\textbf{i}t} & & & & & & & & -1 & & & & & \\
\vdots & \vdots & \vdots & & & & & & & & & \ddots & & & & \\
1 & -\textbf{i}e^{-\textbf{i}t} & \textbf{i}e^{-\textbf{i}t} & & & & & & & & & & -1 & & & \\
\hline
1 & -\textbf{i} & -\textbf{i} & & & & & & & & & & & -1 & & \\
\vdots & \vdots & \vdots & & & & & & & & & & & & \ddots & \\
1 & -\textbf{i} & -\textbf{i} & & & & & & & & & & & & & -1 \\
\end{array}\right]
\]

\normalsize In the figure above the fourth row is  marked as the
one considered. In this row, starting with
$(1,\textbf{i},\textbf{i})$, let $a,c,e$ and $g$ denote the number
of $\textbf{i}$'s, while $b,d,f$ and $h$ denote the number of
$-\textbf{i}$'s in the corresponding ``cells''. Note that by
taking the inner product of the first three rows of $D_N$, one can
calculate how many vertical pairs $(\textbf{i},\textbf{i}),
(\textbf{i},-\textbf{i}), (-\textbf{i},\textbf{i})$ and
$(-\textbf{i},-\textbf{i})$ there can be in rows $(\vu,\vv)$. The following equations are necessary and sufficient conditions for the orthogonality of the first three rows of $D_N^{(1)}(t)$, independently of $t$.

\beql{e1}
b=\frac{N-2}{4}-2-a
\eeq

\beql{e2}
d=\frac{N-2}{4}-c
\eeq

\beql{e3}
f=\frac{N-2}{4}-e
\eeq

\beql{e4}
h=\frac{N-2}{4}-g
\eeq

The number of $\textbf{i}$'s is $\frac{N-2}{2}$ in every row, so we have

\beql{g3}
a+c+e+g=\frac{N-2}{2}-2
\eeq

Since $D_N^{(1)}(0)$ is Hadamard  the fourth row is orthogonal to
$\vu$, prior to modification, and we get

\beql{f1} 2+a-b+c-d-e+f-g+h=0 \eeq

Now put \eqref{e1}-\eqref{e4} into \eqref{f1}, yielding

\beql{g1} a+c=e+g-2 \eeq

Similarly, the fourth row of $D_N^{(1)}(0)$ is orthogonal to
$\vv$,  prior to modification, and we get

\beql{f2}
2+a-b-c+d+e-f-g+h=0
\eeq

Substituting \eqref{e1}-\eqref{e4} into \eqref{f2} implies

\beql{g2}
a+e=c+g-2
\eeq

Finally, use \eqref{g3} in \eqref{g1} and \eqref{g2} to obtain

\beql{}
a+c=a+e \left(=\frac{N-10}{4}\right)
\eeq

This last equation implies that $c=e$, and from \eqref{e2}  and
\eqref{e3} $d=f$ immediately follows. Now it is only a matter of
simple computation, to show that both $\vu$ and $\vv$ are
orthogonal to the chosen row of $D_N^{(1)}(t)$, independently of
the value of $t$.

\textbf{CASE 2:} We need to prove that a row with
$e^{-\textbf{i}t}$-type parameters is orthogonal to both $\vu$ and
$\vv$. Consider a row starting with
$\left(1,\textbf{i}e^{-\textbf{i}t},-\textbf{i}e^{-\textbf{i}t}\right)$
(the case
$\left(1,-\textbf{i}e^{-\textbf{i}t},\textbf{i}e^{-\textbf{i}t}\right)$
can be treated similarly). The columns of $D_N^{(1)}(t)$ can be
permuted so that it takes the form:

\scriptsize

\[
\begin{array}{r}
(\vu)\\
(\vv)\\
\end{array}
\left[\begin{array}{r|rr|r|rrr|rrr|rrr|rrr}
1 & 1 & 1 & 1 & 1 & \hdots & 1 & 1 & 1 & 1 & 1 & 1 & 1 & 1 & 1 & 1\\
\hline
1 & -1 & \textbf{i} & \textbf{i}e^{\textbf{i}t} & \textbf{i} & \hdots & \textbf{i} & \textbf{i}e^{\textbf{i}t} & \hdots & \textbf{i}e^{\textbf{i}t} & -\textbf{i}e^{\textbf{i}t} & \hdots & -\textbf{i}e^{\textbf{i}t} & -\textbf{i} & \hdots & -\textbf{i}\\
1 & \textbf{i} & -1 & -\textbf{i}e^{\textbf{i}t} & \textbf{i} & \hdots & \textbf{i} & -\textbf{i}e^{\textbf{i}t} & \hdots & -\textbf{i}e^{\textbf{i}t} & \textbf{i}e^{\textbf{i}t} & \hdots & \textbf{i}e^{\textbf{i}t} & -\textbf{i} & \hdots & -\textbf{i} \\
\hline
1 & \textbf{i}e^{-\textbf{i}t} & -\textbf{i}e^{-\textbf{i}t} & -1 & $a$ & & $b$ & $c$ & & $d$ & $e$ & & $f$ & $g$ & & $h$\\
\end{array}\right]\]

\normalsize where the fourth row is the one under consideration,
and  $a,b,c,d,e,f,g$ and $h$ have the same meaning as in CASE 1.
Again, we express the orthogonality of the first three rows of $D_N^{(1)}(t)$ as:

\beql{t1}
b=\frac{N-2}{4}-1-a
\eeq

\beql{t2}
d=\frac{N-2}{4}-1-c
\eeq

\beql{t3}
f=\frac{N-2}{4}-e
\eeq

\beql{t4}
h=\frac{N-2}{4}-g
\eeq

And the allowed number of $\textbf{i}$'s is

\beql{t5}
a+c+e+g=\frac{N-2}{2}-1
\eeq

Again, as $\vu$ and $\vv$ are orthogonal to the considered parametrized row for $t=0$, one gets

\beql{t6}
a-b+c-d-e+f-g+h=0
\eeq

and

\beql{t7}
2+a-b-c+d+e-f-g+h=0
\eeq

By substituting \eqref{t1}-\eqref{t4} into \eqref{t6} and \eqref{t7} we get

\beql{t8}
a+c=e+g-1
\eeq

and

\beql{t9}
a+e=c+g-1
\eeq

Again, use \eqref{t5} in \eqref{t8} and \eqref{t9}  to obtain

\beql{}
a+c=a+e\left(=\frac{N-6}{4}\right)
\eeq

This last equation  implies $c=e$ and from \eqref{t2} and
\eqref{t3} $d=f-1$ follows. By applying these identities it is
only a matter of simple computation that the considered
$e^{\textbf{i}t}$-type row is orthogonal to $\vu$ and $\vv$,
independently of $t$.

\textbf{OTHER CASES:} Considering any other pair of rows in
$D_N^{(1)}(t)$ it is trivial to show that they are orthogonal to
each other.  This completes the proof.
\end{proof}

The last theorem allows introduction of one parameter for every
complex  Hadamard matrix arising from the conference matrix
construction. However the following more complex method seems to
be working in general. In some sense this is a natural
generalization of Theorem \ref{b}.

\begin{construction}\label{conj}
Take an arbitrary dephased, symmetric complex Hadamard matrix $D$
arising from the conference matrix construction, of order $N$. Use
Theorem \ref{CO} method, involving a pair of rows (and the corresponding columns), to introduce a free parameter in $D$. Then select another pair of ``suitable'' rows (and the corresponding columns), if possible, in order to use Theorem \ref{CO} again to introduce another parameter. A ``suitable'' pair of rows must satisfy two conditions:
\begin{itemize}
\item[i)] all its vertical pairs of entries are formed (taking into account already existing parameters, if any)
either by identical entries or entries being negative with respect to each other (except for the inevitable $(-1,\ast)$ and $(\ast,-1)$ pairs);
\item[ii)] it has a vertical pair $(\textbf{i},-\textbf{i})$ or $(-\textbf{i},\textbf{i})$, not yet parametrized.
\end{itemize}
If a suitable pair of rows is found, introduce a new parameter in
it (and in the corresponding columns) in the manner analogous to that of Theorem \ref{CO}, i.e. multiplying pairs of opposite entries by $e^{\pm\textbf{i}t}$. Repeat this procedure as long as there exist suitable
pairs of rows.
\end{construction}

The two conditions above seem to be {\it necessary} in the
following sense. Condition i) guarantees that the first row of $D$ and the rows of a newly parametrized pair are all orthogonal to each other, while condition ii) is required to ensure that the
newly introduced parameter does not depend on earlier ones. It is
not clear, however, that they are indeed sufficient, i.e. we do
not have a formal proof that the arising parametric matrices
remain Hadamard. Also, if several suitable pairs of rows exist at
one stage then it is not clear which pair to favour over the
others. The maximal number of parameters that can be introduced in
this way is $\frac{N}{2}-1$ (because the first row definitely does
not have a companion to make a pair with). We used this
construction to obtain the families stemming from $D_{10}$ and
$D_{14}$ below, and the well-known family $D_6^{(1)}(t)$ of \cite{karol} also arises in this way. These examples suggest the following

\begin{conjecture}
Construction \ref{conj} leads to Hadamard matrices after each
step, and for $N\geq 14$ the maximum number, $\frac{N}{2}-1$, of
parameters can be introduced.
\end{conjecture}

\begin{rem}
The construction yields only $\frac{N}{2}-2$ parameters for $D_6$
and $D_{10}$, because condition ii) fails to hold due to the
matrices being ``too small''.
\end{rem}

In the recent catalogue  \cite{karol} only Di\c{t}\u{a}-type
matrices were considered in dimensions $N=10$ and $14$. In view of
Theorem \ref{CO2} and \ref{CO} we can now present new parametric
families of complex Hadamard matrices of these orders. Our first
example is the matrix $D_{10}$ which is constructed from the only
(up to equivalence) conference matrix of order $10$.

\beql{}
D_{10}=\left[\begin{array}{rrrrrrrrrr}
 1 & 1 & 1 & 1 & 1 & 1 & 1 & 1 & 1 & 1 \\
 1 & -1 & -\textbf{i} & -\textbf{i} & -\textbf{i} & -\textbf{i} & \textbf{i} & \textbf{i} & \textbf{i} & \textbf{i} \\
 1 & -\textbf{i} & -1 & \textbf{i} & \textbf{i} & -\textbf{i} & -\textbf{i} & -\textbf{i} & \textbf{i} & \textbf{i} \\
 1 & -\textbf{i} & \textbf{i} & -1 & -\textbf{i} & \textbf{i} & -\textbf{i} & \textbf{i} & -\textbf{i} & \textbf{i} \\
 1 & -\textbf{i} & \textbf{i} & -\textbf{i} & -1 & \textbf{i} & \textbf{i} & -\textbf{i} & \textbf{i} & -\textbf{i} \\
 1 & -\textbf{i} & -\textbf{i} & \textbf{i} & \textbf{i} & -1 & \textbf{i} & \textbf{i} & -\textbf{i} & -\textbf{i} \\
 1 & \textbf{i} & -\textbf{i} & -\textbf{i} & \textbf{i} & \textbf{i} & -1 & -\textbf{i} & -\textbf{i} & \textbf{i} \\
 1 & \textbf{i} & -\textbf{i} & \textbf{i} & -\textbf{i} & \textbf{i} & -\textbf{i} & -1 & \textbf{i} & -\textbf{i} \\
 1 & \textbf{i} & \textbf{i} & -\textbf{i} & \textbf{i} & -\textbf{i} & -\textbf{i} & \textbf{i} & -1 & -\textbf{i} \\
 1 & \textbf{i} & \textbf{i} & \textbf{i} & -\textbf{i} & -\textbf{i} & \textbf{i} & -\textbf{i} & -\textbf{i} & -1
\end{array}\right]
\eeq

We have  already seen that $D_{10}$ is a non-Di\c{t}\u{a}-type
matrix and according to Theorem \ref{CO} it has an affine orbit
stemming from it. Moreover, by Construction \ref{conj} we could
introduce $3$  parameters (we chose the ``suitable'' pairs of rows
by an ad hoc method, as follows: $(2,10),(3,9)$ and $(5,7)$).

\beql{}
D_{10}^{(3)}(a,b,c)=D_{10}\circ \textbf{EXP}\left(\textbf{i}\cdot R_{D_{10}^{(3)}}(a,b,c)\right)
\eeq

where

\beql{}
R_{D_{10}^{(3)}}(a,b,c)=\left[\begin{array}{llllllllll}
 \bullet & \bullet & \bullet & \bullet & \bullet & \bullet & \bullet & \bullet & \bullet & \bullet \\
 \bullet & \bullet & a-b & a & -c & \bullet & -c & a & a-b & \bullet \\
 \bullet & -a+b & \bullet & b & -c & \bullet & -c & b & \bullet & -a+b \\
 \bullet & -a & -b & \bullet & \bullet & \bullet & \bullet & \bullet & -b & -a \\
 \bullet & c & c & \bullet & \bullet & \bullet & \bullet & \bullet & c & c \\
 \bullet & \bullet & \bullet & \bullet & \bullet & \bullet & \bullet & \bullet & \bullet & \bullet \\
 \bullet & c & c & \bullet & \bullet & \bullet & \bullet & \bullet & c & c \\
 \bullet & -a & -b & \bullet & \bullet & \bullet & \bullet & \bullet & -b & -a \\
 \bullet & -a+b & \bullet & b & -c & \bullet & -c & b & \bullet & -a+b \\
 \bullet & \bullet & a-b & a & -c & \bullet & -c & a & a-b & \bullet
\end{array}\right]
\eeq

We checked with a computer that $D_{10}^{(3)}(a,b,c)$ is indeed
Hadamard. The defect (in the sense of \cite{karol}) of $D_{10}$ is
$16$, so we cannot be sure that $D_{10}^{(3)}(a,b,c)$ is maximal
affine (the defect is an upper bound for the dimensionality of a
family stemming from $D_{10}$). It is possible that further
parameters can be introduced.

Now we turn to $N=14$. Our starting point Hadamard matrix,
constructed from the only (up to equivalence) conference matrix of
order $14$,  is the following $D_{14}$.

\beql{}
D_{14}=\left[\begin{array}{rrrrrrrrrrrrrr}
 1 & 1 & 1 & 1 & 1 & 1 & 1 & 1 & 1 & 1 & 1 & 1 & 1 & 1 \\
 1 & -1 & \textbf{i} & -\textbf{i} & \textbf{i} & \textbf{i} & -\textbf{i} & -\textbf{i} & -\textbf{i} & -\textbf{i} & \textbf{i} & \textbf{i} & -\textbf{i} & \textbf{i} \\
 1 & \textbf{i} & -1 & \textbf{i} & -\textbf{i} & \textbf{i} & \textbf{i} & -\textbf{i} & -\textbf{i} & -\textbf{i} & -\textbf{i} & \textbf{i} & \textbf{i} & -\textbf{i} \\
 1 & -\textbf{i} & \textbf{i} & -1 & \textbf{i} & -\textbf{i} & \textbf{i} & \textbf{i} & -\textbf{i} & -\textbf{i} & -\textbf{i} & -\textbf{i} & \textbf{i} & \textbf{i} \\
 1 & \textbf{i} & -\textbf{i} & \textbf{i} & -1 & \textbf{i} & -\textbf{i} & \textbf{i} & \textbf{i} & -\textbf{i} & -\textbf{i} & -\textbf{i} & -\textbf{i} & \textbf{i} \\
 1 & \textbf{i} & \textbf{i} & -\textbf{i} & \textbf{i} & -1 & \textbf{i} & -\textbf{i} & \textbf{i} & \textbf{i} & -\textbf{i} & -\textbf{i} & -\textbf{i} & -\textbf{i} \\
 1 & -\textbf{i} & \textbf{i} & \textbf{i} & -\textbf{i} & \textbf{i} & -1 & \textbf{i} & -\textbf{i} & \textbf{i} & \textbf{i} & -\textbf{i} & -\textbf{i} & -\textbf{i} \\
 1 & -\textbf{i} & -\textbf{i} & \textbf{i} & \textbf{i} & -\textbf{i} & \textbf{i} & -1 & \textbf{i} & -\textbf{i} & \textbf{i} & \textbf{i} & -\textbf{i} & -\textbf{i} \\
 1 & -\textbf{i} & -\textbf{i} & -\textbf{i} & \textbf{i} & \textbf{i} & -\textbf{i} & \textbf{i} & -1 & \textbf{i} & -\textbf{i} & \textbf{i} & \textbf{i} & -\textbf{i} \\
 1 & -\textbf{i} & -\textbf{i} & -\textbf{i} & -\textbf{i} & \textbf{i} & \textbf{i} & -\textbf{i} & \textbf{i} & -1 & \textbf{i} & -\textbf{i} & \textbf{i} & \textbf{i} \\
 1 & \textbf{i} & -\textbf{i} & -\textbf{i} & -\textbf{i} & -\textbf{i} & \textbf{i} & \textbf{i} & -\textbf{i} & \textbf{i} & -1 & \textbf{i} & -\textbf{i} & \textbf{i} \\
 1 & \textbf{i} & \textbf{i} & -\textbf{i} & -\textbf{i} & -\textbf{i} & -\textbf{i} & \textbf{i} & \textbf{i} & -\textbf{i} & \textbf{i} & -1 & \textbf{i} & -\textbf{i} \\
 1 & -\textbf{i} & \textbf{i} & \textbf{i} & -\textbf{i} & -\textbf{i} & -\textbf{i} & -\textbf{i} & \textbf{i} & \textbf{i} & -\textbf{i} & \textbf{i} & -1 & \textbf{i} \\
 1 & \textbf{i} & -\textbf{i} & \textbf{i} & \textbf{i} & -\textbf{i} & -\textbf{i} & -\textbf{i} & -\textbf{i} & \textbf{i} & \textbf{i} & -\textbf{i} & \textbf{i} & -1
\end{array}\right]
\eeq

Again, this is a non-Di\c{t}\u{a}-type matrix,  and a
$6$-parameter affine family stems from it (which we constructed
with the aid of Construction \ref{conj}; the considered ``suitable''
pairs of rows were $(2,3),(4,5),(6,9),(7,13),(8,12)$ and
$(11,14)$). The defect of the matrix is $36$ so it might be
possible to introduce further parameters. We do not claim that
\emph{all} the matrices contained in the family stemming from
$D_{14}$ are non-Di\c{t}\u{a}-type, but it is obviously true in a
small neighborhood of it.

\beql{}
D_{14}^{(6)}(a,b,c,d,e,f)=D_{14}\circ \textbf{EXP}\left(\textbf{i}\cdot R_{D_{14}^{(6)}}(a,b,c,d,e,f)\right)
\eeq

where

\beql{}
R_{D_{14}^{(6)}}(a,b,c,d,e,f)=
\eeq
\scriptsize
\[
\left[
\begin{array}{cccccccccccccc}
 \bullet & \bullet & \bullet & \bullet & \bullet & \bullet & \bullet & \bullet & \bullet & \bullet & \bullet & \bullet & \bullet & \bullet \\
 \bullet & \bullet & \bullet & a-b & a-b & -c & a & -e & -c & \bullet & a & -e & a & a \\
 \bullet & \bullet & \bullet & a-b & a-b & -c & a & -e & -c & \bullet & a & -e & a & a \\
 \bullet & b-a & b-a & \bullet & \bullet & b & b & -e & b & \bullet & -f & -e & b & -f \\
 \bullet & b-a & b-a & \bullet & \bullet & b & b & -e & b & \bullet & -f & -e & b & -f \\
 \bullet & c & c & -b & -b & \bullet & c-d & c & \bullet & \bullet & \bullet & c & c-d & \bullet \\
 \bullet & -a & -a & -b & -b & d-c & \bullet & d-e & d-c & \bullet & d-f & d-e & \bullet & d-f \\
 \bullet & e & e & e & e & -c & e-d & \bullet & -c & \bullet & -f & \bullet & e-d & -f \\
 \bullet & c & c & -b & -b & \bullet & c-d & c & \bullet & \bullet & \bullet & c & c-d & \bullet \\
 \bullet & \bullet & \bullet & \bullet & \bullet & \bullet & \bullet & \bullet & \bullet & \bullet & \bullet & \bullet & \bullet & \bullet \\
 \bullet & -a & -a & f & f & \bullet & f-d & f & \bullet & \bullet & \bullet & f & f-d & \bullet \\
 \bullet & e & e & e & e & -c & e-d & \bullet & -c & \bullet & -f & \bullet & e-d & -f \\
 \bullet & -a & -a & -b & -b & d-c & \bullet & d-e & d-c & \bullet & d-f & d-e & \bullet & d-f \\
 \bullet & -a & -a & f & f & \bullet & f-d & f & \bullet & \bullet & \bullet & f & f-d & \bullet
\end{array}
\right]
\]
\normalsize

To summarize the cases $N=10, 14$ we conclude that

\begin{corollary}
The families $D_{10}^{(3)}(a,b,c)$ and $D_{14}^{(6)}(a,b,c,d,e,f)$
are locally inequivalent to the families contained in
\cite{karol}.
\end{corollary}

\begin{rem}
Note that $D_{10}$ and $D_{14}$ are unique in the sense that
according to \cite{Cr} the number of  inequivalent symmetric
conference matrices is $1$ for orders $N=2,6,10,14$ and $18$,
while already for order $N=26$ there exist $4$ inequivalent
symmetric conference matrices. This implies that in higher
dimensions it may be possible to construct locally inequivalent families stemming from inequivalent starting point matrices. Recall that there is no conference matrix of order $22$ and $34$ due to Raghavarao's theorem \cite{R}.
\end{rem}

Let us summarize our results. In this paper we have described two
general constructions of parametric families of complex
Hadamard matrices. We have presented new matrices of order $10,12$
and $14$, thus we have supplemented the recent catalogue of complex
Hadamard matrices of small orders \cite{karol}. We pointed out
that certain real Hadamard matrices cannot be constructed using
Di\c{t}\u{a}'s formula, so in order to find all inequivalent
complex Hadamard matrices of a given order one should look for and
resort to other construction methods.

It would be interesting  to see whether the hereby presented
families can be extended with further parameters. It also remains to
be checked whether Construction \ref{conj} leads indeed to
parametric families of complex Hadamard matrices in general.

\section*{Acknowledgement}
The author thanks M\'at\'e Matolcsi for many insightful comments.
He is also greatly indebted to the referee whose suggestions have
substantially improved the presentation of the paper.

\end{document}